%&amstex
\input amstex
\input amsppt.sty
\magnification=\magstep1
\hsize=30truecc
\vsize=22.2truecm
\baselineskip=16truept
%\NoBlackBoxes
\TagsOnRight
\nologo
\pageno=1
\topmatter
\def\N{\Bbb N}
\def\Z{\Bbb Z}
\def\Q{\Bbb Q}

\def\C{\Bbb C}
\def\l{\left}
\def\r{\right}
\def\b{\bigg}

\def\({\b(}
\def\[{\b[}
\def\){\b)}
\def\]{\b]}

\def\t{\text}
\def\f{\frac}
\def\mo{\roman{mod}}
\def\ord{\roman{ord}}

\def\bi{\binom}
\def\eq{\equiv}

\def\ls{\leqslant}
\def\gs{\geqslant}
\def\al{\alpha}

\def\Proof{\noindent{\it Proof}}

\hbox {Preprint (August 23, 2006), {\tt arXiv:math.NT/0608560}.}
\bigskip
\title Extensions of Wilson's Lemma and the Ax-Katz Theorem\endtitle
\author  Zhi-Wei Sun\endauthor
\affil Department of Mathematics, Nanjing University
\\ Nanjing 210093, People's Republic of China
\\zwsun\@nju.edu.cn
\\ {\tt http://pweb.nju.edu.cn/zwsun}
\medskip
\endaffil
\thanks 2000 {\it Mathematics Subject Classification}.\,Primary 11T06;
Secondary 05A10, 11A07, 11S05, 41A10.\newline\indent
Supported by the National Science Fund
for Distinguished Young Scholars (No. 10425103) in China.
\endthanks
\abstract A classical result of A. Fleck states that
if $p$ is a prime, and $n>0$ and $r$ are integers, then
$$\sum_{k\eq r\,(\mo\ p)}\bi nk(-1)^k\eq0\ \l(\mo\ p^{\lfloor(n-1)/(p-1)\rfloor}\r).$$
Recently R. M. Wilson used Fleck's congruence and Weisman's extension
to present a useful lemma on polynomials modulo prime powers, and applied this lemma
to reprove the Ax-Katz theorem on solutions of congruences modulo $p$
and deduce various results on codewords in $p$-ary linear codes
with weights. In light of the recent generalizations of Fleck's congruence
given by D. Wan, and D. M. Davis and Z. W. Sun, we obtain new extensions
of Wilson's lemma and the Ax-Katz theorem.
\endabstract
\endtopmatter

\document

\heading{1. Introduction}\endheading

Let $p$ be a prime, and let $n\in\N=\{0,1,2,\ldots\}$ and $r\in\Z$.
In 1913 A. Fleck (cf. [D, p.\,274]) proved that
$$\ord_p\(\sum_{k\eq r\,(\mo\ p)}\bi nk(-1)^k\)\gs\l\lfloor\f{n-1}{p-1}\r\rfloor,\tag1.1$$
where $\lfloor\cdot\rfloor$ is the well-known floor function, and the $p$-adic order
$\ord_p(\al)$ of a $p$-adic number $\al$
is given by $\sup\{a\in\Z:\, \al/p^a\in\Z_p\}$.
(As usual $\Z_p$ denotes the ring of $p$-adic integers in the $p$-adic field $\Q_p$.)

Let $a\in\Z^+=\{1,2,3,\ldots\}$. In 1977, motivated by
his study of $p$-adically continuous functions and unaware of Fleck's earlier result,
C. S. Weisman [We] extended Fleck's inequality as follows:
$$\ord_p\(\sum_{k\eq r\,(\mo\ p^a)}\bi nk(-1)^k\)\gs\l\lfloor\f{n-p^{a-1}}{\varphi(p^a)}\r\rfloor,\tag1.2$$
where $\varphi$ is Euler's totient function.

For a function $f$ from the complex field $\C$ to $\C$, let
$\Delta^0 f(x)=f(x)$,
$\Delta f(x)=f(x+1)-f(x)$ and $\Delta^n f(x)=\Delta\Delta^{n-1}f(x)$
for $n=2,3,\ldots$.
Now we recall a classical interpolation formula due to I. Newton and J. Gregory.

\proclaim{Newton-Gregory Interpolation Formula} Given a function $f:\C\to\C$, for any $d\in\N$ we have
$$f(x)=\sum_{n=0}^d c_n\bi xn+R_d(x),$$
where
$$c_n=\Delta^nf(0)=\sum_{k=0}^n\bi nk(-1)^{n-k}f(k)$$
and
$$R_d(x)=\left|\matrix
1&0&\cdots &0&f(0)\\
1&1^1&\cdots &1^d&f(1)\\
\cdots&\cdots&\cdots&\cdots&\cdots\\
1&d^1&\cdots&d^d&f(d)
\\1&x^1&\cdots&x^d&f(x)
\endmatrix\right|
\bigg/\left|\matrix
1^1&1^2&\cdots&1^d\\2^1&2^2&\cdots&2^d\\
\cdots&\cdots&\cdots&\cdots\\
d^1&d^2&\cdots&d^d
\endmatrix\right|.$$
$($Note that $R_d(x)=0$ if $f$ is a polynomial with $\deg f\ls d)$.
\endproclaim

In 2006 R. M. Wilson [Wi] rediscovered Weisman's (1.2) in the case $n\eq p^{a-1}\ (\mo\ \varphi(p^a))$,
and used it to obtain the following lemma (similar to the Newton-Gregory interpolation formula)
and give many applications.

\proclaim{Wilson's Lemma} Let $p$ be a prime, and let $a,b\in\Z^+$.
Let $f$ be an integer-valued function on the integers that is periodic modulo $p^a$. Then
there exists a polynomial
$$w(x)=c_0+c_1x+c_2\bi x2+\cdots+c_d\bi xd\ \ (c_0,c_1,\ldots,c_d\in\Z)$$
 of degree smaller than $b\varphi(p^a)+p^{a-1}$ such that
$$\ord_p(c_n)\gs\l\lfloor\f{n-p^{a-1}}{\varphi(p^a)}\r\rfloor\ \ \t{for all}\ n=0,\ldots,d,$$
 and $w(x)\eq f(x)\ (\mo\ p^b)$ for
 all $x\in\Z$.
 \endproclaim

In this paper, for a prime $p$ we let $\overline \Q_p$ be the algebraic closure of the field $\Q_p$
and let $\overline\Z_p$ be the ring of $p$-adic algebraic integers in $\overline\Q_p$.
For $m,n\in\N$ we use $[m,n]$ to denote the set $\{x\in\Z:\, m\ls x\ls n\}$.

In view of the recent generalizations of Fleck's and Weisman's results (cf. [S], [W06], [SW], [DS] and [SD]),
we are able to present the following further extension of Wilson's Lemma.

\proclaim{Theorem 1.1} Let $p$ be a prime, and let $a\in\N$ and $b\in\Z^+$.
Let $f(x)\in\overline\Q_p[x]$ with $\deg f\ls l\in\N$ and $f(m)\in\overline\Z_p$ for all $m\in\Z$,
and let $g$ be a function from $[0,p^a-1]$ to $\overline\Z_p$. Let
$d\in\N$ be the maximal integer with $M_d<b$, where $M_d$ denotes
$$\max\l\{\l\lfloor\f{d-lp^a-p^{a-1}}{\varphi(p^a)}\r\rfloor,
\ord_p\(\l\lfloor\f{d}{p^{a-1}}\r\rfloor!\)-\ord_p(l!)-\min\l\{l,\l\lfloor\f d{p^a}\r\rfloor\r\}\r\}.$$
Then there exists
a polynomial
$$P(x)=\sum_{n=0}^dc_n\bi xn\quad (c_0,\ldots,c_d\in\overline\Z_p)\tag1.3$$
with $\ord_p(c_n)\gs M_n$ for all $n=0,\ldots,d$, such that
$$P(p^aq+r)\eq f(q)g(r)\ (\mo\ p^b)\ \ \t{for all}\ q\in\Z\ \t{and}\ r\in[0,p^a-1].\tag1.4$$
\endproclaim

The following celebrated theorem (cf. C. Chevally [C],
E. Warning [Wa] and Theorem 2.6 of M. B. Nathanson [N, pp.\,50--51])
is well known and quite useful.

\proclaim{Chevalley-Warning Theorem} Let
$f_1(x,\ldots,x_n),\ldots,f_m(x_1,\ldots,x_n)$
be polynomials over a finite field $F$ of characteristic $p$
with $\deg f_1+\cdots+\deg f_m<n$. Then
the number of solutions to the system of equations
$$f_1(x_1,\ldots,x_n)=0,\ \ldots,\ f_m(x_1,\ldots,x_n)=0\tag1.5$$
over $F^n$ is a multiple of $p$.
\endproclaim

Here is a further refinement of the Chevalley-Warning theorem due to J. Ax [A] in the case $m=1$,
and N. Katz [K] in the general case.
\proclaim{Ax-Katz Theorem} Let $F_q$ be the finite field with $q=p^a$ elements
where $p$ is a prime and $a\in\Z^+$. Let
$f_1(x,\ldots,x_n),\ldots,f_m(x_1,\ldots,x_n)$
be nonzero polynomials over $F_q$ with degrees $d_1\gs\ldots\gs d_m$ respectively.
Then, for any positive integer $b$ satisfying $n>(b-1)d_1+(d_1+\cdots+d_m)$,
$q^b$ divides the number of solutions to the system $(1.5)$ over $F^n$.
\endproclaim

D. Wan [W89, W95] gave a new proof of the Ax-Katz theorem via the Stickelberger theorem.
In 2005 X.-D. Hou [H] reduced the Ax-Katz theorem to the Ax theorem on a single polynomial equation.
In 2006 Wilson [Wi] reproved the Ax-Katz theorem for prime fields by using Wilson's Lemma.

With help of Theorem 1.1, we establish the following theorem.
\proclaim{Theorem 1.2} Let $p$ be a prime, and let $F_1(x),\ldots,F_m(x)\in\overline\Q_p[x]$
with $\deg F_k\ls l_k\in\N$ and $F_k(a)\in\overline\Z_p$ for all $a\in\Z$.
Let $a_1,\ldots,a_m\in\N$, and let $f_1(x_1,\ldots,x_n),\ldots,f_m(x_1,\ldots,x_n)$
be nonzero polynomials with integer coefficients.
Assume that
$d_1\varphi(p^{a_1})=\max_{1\ls k\ls m}d_k\varphi(p^{a_k})$
where $d_k=\deg f_k$ for $k=1,\ldots,m$.
Let $b\in\Z^+$ and suppose that
$$n>(b-1)\max\l\{\f{d_1\varphi(p^{a_1})}{p-1},1\r\}+\f1{p-1}\sum_{k=1}^m((l_k+1)p^{a_k}-[\![a_k\not=0]\!])d_k,\tag1.6$$
where $[\![a_k\not=0]\!]$ takes $1$ or $0$ according as $a_k\not=0$ or not.
Then
$$\sum\Sb x_1,\ldots,x_n\in[0,p-1]
\\p^{a_k}\mid f_k(x_1,\ldots,x_n)\ \t{for all}\ k\in[1,m]\endSb
\prod_{k=1}^mF_k\l(\f{f_k(x_1,\ldots,x_n)}{p^{a_k}}\r)\eq0\ (\mo\ p^b).\tag1.7$$
\endproclaim

In the case $F_1(x)=\cdots=F_m(x)=1$, Theorem 1.2 yields an extension
of the Ax-Katz theorem for prime fields. In 1995 O. Moreno and C. J. Moreno [MM]
introduced a method to reduce the general case of the Ax-Katz theorem
to the prime field case.

\proclaim{Corollary 1.1} Let $f_1(x_1,\ldots,x_n),\ldots,
f_m(x_1,\ldots,x_n)$ be nonzero polynomials with integer coefficients
having degrees $d_1\gs\cdots\gs d_m$ respectively.
If $p$ is a prime, $a,b\in\Z^+$, $l_1,\ldots,l_m\in\N$ and
$$n>(b-1)d_1p^{a-1}+\f{p^a-1}{p-1}\sum_{k=1}^md_k+\f{p^a}{p-1}\sum_{k=1}^ml_kd_k,\tag1.8$$
then we have
$$\sum\Sb x_1,\ldots,x_n\in[0,p-1]\\p^a\mid f_k(x_1,\ldots,x_n)\ \t{for all}\ k\in[1,m]\endSb
\prod_{k=1}^m\bi{f_k(x_1,\ldots,x_n)/p^a}{l_k}\eq0\ (\mo\ p^b).\tag1.9$$
\endproclaim
\Proof. Just apply Theorem 1.2 with $a_k=a$ and $F_k(x)=\bi x{l_k}$ for $k=1,\ldots,m$. \qed

\medskip
Let $q=p^a$ where $p$ is a prime and $a\in\Z^+$, and let
$\zeta_{q-1}\in\overline\Z_p$ be a primitive $(q-1)$-th roots of unity.
It is well known that $\Z_p[\zeta_{q-1}]/(p)$ is a finite field of $q$ elements.
The finite field $\Bbb F_q$ of $q=p^a$ elements
is an extension of the prime field $\Bbb F_p$ with $[\Bbb F_q:\Bbb F_p]=a$.
Thus $\Bbb F_q$ is isomorphic to $\Bbb F_p^a$ and the Chevalley-Warning theorem
can be reduced to the prime field case. Corollary 1.1 in the case
$a=b=1$ and $l_1=\ldots=l_m=0$ yields the Chevalley-Warning theorem for $\Bbb F_p=\Z/p\Z$
and hence the general case of the Chevalley-Warning theorem.

\heading{2. Proofs of Theorems 1.1 and 1.2}\endheading

\proclaim{Lemma 2.1} Let $p$ be a prime, and let $f(x)\in\overline\Q_p[x]$ with $\deg f\ls l\in\N$
and $f(m)\in\overline\Z_p$ for all $m\in\Z$. For any $a,n\in\N$ and $r\in\Z$, we have
$$\ord_p\(\sum_{k\eq r\,(\mo\ p^a)}\bi nk(-1)^kf\l(\f{k-r}{p^a}\r)\)
\gs\l\lfloor\f{n-lp^a-p^{a-1}}{\varphi(p^a)}\r\rfloor\tag 2.1$$
and
$$\aligned&\ord_p\(\sum_{k\eq r\,(\mo\ p^a)}\bi nk(-1)^kf\l(\f{k-r}{p^a}\r)\)
\\\gs&\ord_p\(\l\lfloor\f{n}{p^{a-1}}\r\rfloor!\)
-\ord_p(l!)-\min\l\{l,\l\lfloor\f n{p^a}\r\rfloor\r\}.
\endaligned\tag 2.2$$
\endproclaim
\Proof. Let $c_j=\sum_{i=0}^j\bi ji(-1)^{j-i}f(i)\in\overline \Z_p$ for $j=0,\ldots,l$.
As $\deg f\ls l$ and $f(x)-\sum_{j=0}^lc_j\bi xj$ vanishes at $0,\ldots,l$, we have
$f(x)=\sum_{j=0}^lc_j\bi xj$.
So it suffices to consider the case $f(x)=\bi xl$ only.

If $a\in\Z^+$ then
$$O:=\ord_p\(\sum_{k\eq r\,(\mo\ p^a)}\bi nk(-1)^k\bi{(k-r)/p^a}l\)
\gs\l\lfloor\f{n-lp^a-p^{a-1}}{\varphi(p^a)}\r\rfloor$$
by D. Wan [W06, Theorem 1.3] (see also [SW] for a combinatorial proof).
This is also true in the case $a=0$, since
$$\sum_{k=0}^n\bi nk(-1)^k\bi{k-r}l=[\![l\gs n]\!](-1)^n\bi{-r}{l-n}$$
 by a known identity (cf. [GKP, (5.24)]).

As $l!\bi xl\in\Z[x]$, by [DS, Theorem 1.5] we have
$$O+\ord_p(l!)\gs\ord_p\(\l\lfloor\f{n}{p^{a}}\r\rfloor!\)=\sum_{s=a+1}^\infty\l\lfloor\f n{p^{s}}\r\rfloor
=\ord_p\(\l\lfloor\f{n}{p^{a-1}}\r\rfloor!\)-\l\lfloor\f n{p^a}\r\rfloor.$$
By [SD, Theorem 1.2], we also have
$$O\gs\ord_p\(\l\lfloor\f{n}{p^{a-1}}\r\rfloor!\)-l-\ord_p(l!).$$

Combining the above we obtain both (2.1) and (2.2). \qed

\medskip
\noindent{\it Proof of Theorem 1.1}. Let $F(x)=f(\lfloor x/p^a\rfloor)g(\{x\}_{p^a})$  for $x\in\Z$,
where $\{x\}_{p^a}$ denotes the least nonnegative residue of $x$ modulo $p^a$. For
$$\align c_n:=&\sum_{k=0}^n\bi nk(-1)^{n-k}F(k)
\\=&(-1)^n\sum_{r=0}^{p^a-1}g(r)\sum_{k\eq r\,(\mo\ p^a)}\bi nk(-1)^kf\l(\f{k-r}{p^a}\r),
\endalign$$
we have $\ord_p(c_n)\gs M_n$ by Lemma 2.1.
If $n>d$, then $\ord_p(c_n)\gs M_n\gs b$.
Set $P(x)=\sum_{n=0}^dc_n\bi xn$. Then, for each $m\in\N$ we have
$$\align F(m)=&\sum_{k=0}^m\bi mkF(k)(1-1)^{m-k}=\sum_{k=0}^m\bi mkF(k)\sum_{n=k}^m\bi{m-k}{n-k}(-1)^{n-k}
\\=&\sum_{n=0}^m\bi mn\sum_{k=0}^n\bi nk(-1)^{n-k}F(k)=\sum_{n\in\N}\bi mnc_n
\\\eq&\sum_{n=0}^d\bi mnc_n=P(m)\ (\mo\ p^b).
\endalign$$
Therefore $P(p^aq+r)\eq F(p^aq+r)=f(q)g(r)\ (\mo\ p^b)$
for all $q\in\N$ and $r\in[0,p^a-1]$.

Choose $N\in\N$ such that $N-b\gs\ord_p(k)$ for all $k\in[1,\max\{d,l\}]$.
For any $x\in\Z$ and $n\in[0,\max\{d,l\}]$, by the Chu-Vandermonde convolution identity
(cf. [GKP, (5.27)]) we have
$$\align&\bi{x+p^N}n=\sum_{k=0}^n\bi{p^N}k\bi x{n-k}
\\=&\bi xn+\sum_{0<k\ls n}\f {p^N}{k}\bi{p^N-1}{k-1}\bi x{n-k}\eq\bi xn\ (\mo\ p^b).
\endalign$$
Therefore $P(x+p^N)\eq P(x)\ (\mo\ p^b)$ and $f(x+p^N)\eq f(x)\ (\mo\ p^b)$
for all $x\in\Z$. For $m=-p^aq+r$ with $q\in\Z^+$ and $r\in[0,p^a-1]$,
clearly $m+p^{a+q+N}\gs0$ and hence
$$\align P(m)\eq &P(m+p^{a+q+N})\eq F(m+p^{a+q+N})
\\\eq& f\l(\l\lfloor\f m{p^a}\r\rfloor+p^{q+N}\r)g(\{m\}_{p^a})
\\\eq& f\l(\l\lfloor\f m{p^a}\r\rfloor\r)g(\{m\}_{p^a})=F(m)\ (\mo\ p^b).
\endalign$$

By the above, we do have $P(p^aq+r)\eq F(p^aq+r)=f(q)g(r)\ (\mo\ p^b)$
for all $q\in\Z$ and $r\in[0,p^a-1]$. \qed

\proclaim{Lemma 2.2} Let $p$ be a prime, and let
$$F(x_1,\ldots,x_n)=\bi{f_1(x_1,\ldots,x_n)}{j_1}\cdots\bi{f_m(x_1,\ldots,x_n)}{j_m},$$
where $j_k\in\N$ and
$f_k(x_1,\ldots,x_n)\in\overline\Z_p[x_1,\ldots,x_n]$ for $k=1,\ldots,m$.
If the total degree of $F(x_1,\ldots,x_n)$ is smaller than $(n-c+1)(p-1)$
for some $c\in\N$, then
$$\sum_{x_1,\ldots,x_n\in[0,p-1]}F(x_1,\ldots,x_n)\eq0\ (\mo\ p^c).$$
\endproclaim
\Proof. See Lemma 4 of Wilson [Wi] and its proof. \qed

\medskip
\noindent{\it Proof of Theorem 1.2}. Given $k\in[1,m]$, by Theorem 1.1
there is a polynomial
$$P_k(x)=\sum_{j=0}^{n_k}c^{(k)}_j\bi xj\quad(c^{(k)}_1,\ldots,c^{(k)}_{n_k}\in\overline\Z_p)$$
such that
$$\ord_p(c^{(k)}_{j})\gs\l\lfloor\f{j-l_kp^{a_k}-p^{a_k-1}}{\varphi(p^{a_k})}\r\rfloor$$
for all $j=0,\ldots,n_k$, and
$$P_k(x)\eq[\![p^{a_k}\mid x]\!]F_k\l(\f x{p^{a_k}}\r)\ (\mo\ p^b)\quad\t{for all}\ x\in\Z.$$
Therefore
$$\align&\sum\Sb x_1,\ldots,x_n\in[0,p-1]\\
p^{a_k}\mid f_k(x_1,\ldots,x_n)\ \t{for all}\ k\in[1,m]\endSb\prod_{k=1}^mF_k\l(\f{f_k(x_1,\ldots,x_n)}{p^{a_k}}\r)
\\\eq&\sum_{x_1,\ldots,x_n\in[0,p-1]}\prod_{k=1}^mP_k(f_k(x_1,\ldots,x_n))
\\\eq&\sum_{j_1=0}^{n_1}c^{(1)}_{j_1}\cdots\sum_{j_m=0}^{n_m}c^{(m)}_{j_m}S(j_1,\ldots,j_m)\ (\mo\ p^b),
\endalign$$
where
$$S(j_1,\ldots,j_m)=\sum_{x_1,\ldots,x_n\in[0,p-1]}\prod_{k=1}^m\bi{f_k(x_1,\ldots,x_n)}{j_k}.$$

Fix $j_1\in[0,n_1],\ldots,j_m\in[0,n_m]$, and let
$$\al_k=\max\l\{\l\lfloor\f{j_k-l_kp^{a_k}-p^{a_k-1}}{\varphi(p^{a_k})}\r\rfloor,0\r\}\quad\t{for}\ k=1,\ldots,m.$$
Then
$$\ord_p\l(c^{(1)}_{j_1}\cdots c_{j_m}^{(m)}\r)=\sum_{k=1}^m\ord_p\l(c^{(k)}_{j_k}\r)\gs\sum_{k=1}^m\al_k.$$
So it suffices to show that $\ord_p(S(j_1,\ldots,j_m))\gs c=b-\sum_{k=1}^m\al_k$.

Assume that $c>0$. By the definition of $\al_k$, $j_k-l_kp^{a_k}-p^{a_k-1}<(\al_k+1)\varphi(p^{a_k})$
and hence
$$j_k\ls l_kp^{a_k}+(\al_k+1)\varphi(p^{a_k})+[\![a_k\not=0]\!](p^{a_k-1}-1).$$
Thus
$$\align\sum_{k=1}^mj_kd_k\ls&\sum_{k=1}^m\l(l_kp^{a_k}+[\![a_k\not=0]\!](p^{a_k-1}-1)+(\al_k+1)\varphi(p^{a_k})\r)d_k
\\=&\sum_{k=1}^m\l(l_kp^{a_k}+p^{a_k}-[\![a_k\not=0]\!]+\al_k\varphi(p^{a_k})\r)d_k
\\\ls&\sum_{k=1}^m\l((l_k+1)p^{a_k}-[\![a_k\not=0]\!]\r)d_k+\varphi(p^{a_1})d_1\sum_{k=1}^m\al_k
\endalign$$
and hence
$$\align \sum_{k=1}^mj_kd_k<&n(p-1)-(b-1)\max\l\{d_1\varphi(p^{a_1}),p-1\r\}+(b-c)d_1\varphi(p^{a_1})
\\\ls&n(p-1)-(c-1)\max\l\{d_1\varphi(p^{a_1}),p-1\r\}.
\endalign$$
Therefore
$$\deg\prod_{k=1}^m\bi{f_k(x_1,\ldots,x_n)}{j_k}\ls\sum_{k=1}^mj_kd_k<(p-1)(n-c+1)$$
and hence $S(j_1,\ldots,j_m)\eq0\ (\mo\ p^c)$ by Lemma 2.2.
This concludes the proof. \qed


\bigskip

\widestnumber\key{GKP}

\Refs

\ref\key A\by J. Ax\paper Zeroes of polynomials over finite fields
\jour Amer. J. Math.\vol 86\yr 1964\pages 255--261\endref

\ref\key C\by C. Chevalley\paper D\'emonstration d'une hypoth\`ese de M. Artin
\jour Abh. Math. Sem. Hamburg\vol 11\yr 1936\pages 73--75\endref

\ref\key DS\by D. M. Davis and Z. W. Sun\paper A number-theoretic
approach to homotopy exponents of SU$(n)$ \jour J. Pure Appl.
Algebra, in press. Available from the website {\tt
http://arxiv.org/abs/math.AT/0508083}\endref

\ref\key D\by  L. E. Dickson\book
History of the Theory of Numbers, {\rm Vol. I}
\publ AMS Chelsea Publ., 1999\endref

\ref\key GKP\by R. Graham, D. E. Knuth and O. Patashnik
 \book Concrete Mathematics
 \publ Addison-Wesley, New York\yr 1989\endref

\ref\key H\by X.-D. Hou\paper A note on the proof of a theorem of Katz
\jour Finite Fields Appl.\vol 11\yr 2005\pages 316--319\endref

\ref\key K\by N. Katz\paper On a theorem of Ax\jour Amer. J. Math.\vol 93\yr 1971\pages 485--499\endref

\ref\key MM\by O. Moreno and C. J. Merono
\paper Improvements of the Chevalley-Warning and the Ax-Katz theorem
\jour Amer. J. Math.\vol 117\yr 1995\pages 241--244\endref

\ref\key N\by M. B. Nathanson\book Additive Number Theory: Inverse Problems and the
Geometry of Sumsets {\rm (Graduate texts in mathematics; 165)}
\publ Springer-Verlag, New York\yr 1996\endref

\ref\key S\by Z. W. Sun\paper Polynomial extension of Fleck's
congruence\jour Acta Arith. \vol 122\yr 2006\pages 91--100\endref

\ref\key SD\by Z. W. Sun and D. M. Davis\paper Combinatorial congruences
modulo prime powers \jour Trans. Amer. Math. Soc.,
in press, {\tt http://arxiv.org/abs/math.NT/0508087}\endref

\ref\key SW\by Z. W. Sun and D. Wan\paper Lucas type congruences
for cyclotomic $\psi$-coefficients\jour preprint, 2005. On-line version:
 {\tt http://arxiv.org/abs/math.NT/0512012}\endref

\ref\key W89\by D. Wan\paper An elementary proof of a theorem of Katz
\jour Amer. J. Math.\vol 111\yr 1989\pages 1--8\endref

\ref\key W95\by D. Wan\paper A Chevalley-Warning approach to $p$-adic estimates of character sums
\jour Proc. Amer. Math. Soc.\vol 123\yr 1995\pages 45--54\endref

\ref\key W06\by D. Wan\paper Combinatorial congruences and $\psi$-operators
\jour Finite Fields Appl.\finalinfo in press,
{\tt http://arxiv.org/abs/math.NT/0603462}\endref

\ref\key Wa\by E. Warning\paper Bemerkung zur vorstehenden Arbeit von Herrn Chevalley
\jour Abh. Math. Sem. Hamburg\vol 11\yr 1936\pages 76--83\endref

\ref\key We\by C. S. Weisman\paper Some congruences for binomial coefficients
\jour Michigan Math. J.\vol 24\yr 1977\pages 141--151\endref

\ref\key Wi\by R. M. Wilson\paper A lemma on polynomials modulo $p^m$
and applications to coding theory\jour Discrete Math.\finalinfo
in press\endref

\endRefs
\enddocument